\documentclass[preprint,review,10pt]{elsarticle}
\usepackage{lineno,hyperref}
\usepackage{graphicx}
\usepackage{subfigure}
\usepackage{amsmath}
\usepackage{geometry}
\usepackage{booktabs}
\usepackage{threeparttable}
\usepackage{textcomp}
\usepackage{multirow}
\usepackage{color,xcolor}
\usepackage{amssymb}
\usepackage{changepage}
\usepackage[latin1]{inputenc}

\modulolinenumbers[50]
\journal{ArXiv Submission Report}

\begin{document}

\begin{frontmatter}

\title{Fixed Point Theorems for Upper Semicontinuous Set-valued Mappings in p-Vector Spaces}


\author{George Xianzhi Yuan}
\address{Chengdu University, Chengdu 610601, China}
\address{College of Mathematics, Sichuan University, Chengdu 610065, China}
\address{Sun Yat-Sen University, Guangzhou 510275, China, and }
\address{East China University of Science and Technology, Shanghai 200237, China}
\address{\bf george\_yuan99@yahoo.com}
\address{This is the revised version of the paper with the same title (arXiv:2303.07177v1 [math.FA], December 13, 2022)}

\begin{abstract}

The goal of this paper is to establish a general fixed point theorem for compact single-valued continuous mapping in Hausdorff $p$-vector spaces, and the fixed point theorem for upper semicontinuous set-valued mappings in Hausdorff locally $p$-convex for $p\in (0, 1]$. These new results provide an answer to Schauder conjecture in the affirmative  under the setting of general $p$-vector spaces for compact single-valued continuous, and also give the fixed point theorems for upper semicontinuous set-valued mappings defined on $s$-convex subsets in Hausdorff locally $p$-convex spaces, which would be fundamental for nonlinear functional analysis in mathematics, where $s, p \in (0. 1]$.

\end{abstract}

\vskip.1in

\end{frontmatter}

\newpage
\linenumbers
\section{Introduction}
\label{intro}

\vskip.1in
It is known that the class of $p$-seminorm spaces $(0 < p \leq 1)$ is an important generalization of usual normed spaces with  rich topological and geometrical structures, and related study has received a lot of attention (e.g., see Alghamdi et al.\cite{alghamdiporegan}, Balachandran \cite{balachandra}, Bayoumi \cite{bayoumi},  Ennassik and Taoudi \cite{ennassik2021}, Ennassik et al.\cite{ennassik}, Gholizadeh et al.\cite{gholizadeh}, Granas and Dugundji \cite{granas}, Jarchow \cite{jarchow}, Kalton \cite{kalton1977}-\cite{kalton1977b}, Kalton et al.\cite{kalton1984}, Park \cite{park100}, Rolewicz \cite{rolewicz}, Xiao and Lu \cite{xiaolu}, Xiao and Zhu \cite{xiaozhu2011}, Yuan \cite{yuan1998}-\cite{yuan2022c}, and many others). However, to the best of our knowledge, the corresponding basic tools and associated results in the category of nonlinear functional analysis have not been well developed, thus the goal of this paper is to develop some
fundamental fixed point theorems for both single-valued and upper semicontinuous set-valued mappings under the framework of Hausdorff and $p$-vector spaces and locally $p$-convex spaces. which include the topological vector spaces and locally convex spaces as special classes, respectively.

We all know that Schauder's fixed point theorem \cite{schauder} in normed spaces is one of the most powerful
tools in dealing with nonlinear problems in analysis. Most notably, it has played a major role in the
development of fixed point theory and related nonlinear analysis and mathematical theory of partial and
differential equations and others. A generalization of Schauder's theorem from normed space to general
topological vector spaces is an old conjecture in fixed point theory which is explained by the Problem 54 of the book $``$The Scottish Book" by Mauldin \cite{Mauldin} as stated as Schauder's conjecture: $``${\bf Every nonempty compact convex set in a topological vector space has the fixed point property};  or in its analytic way, {\bf does a continuous function defined on a compact convex subset of a topological vector space to itself have a fixed point?}"

By following the brief discussion by Ennassik and Taoudi \cite{ennassik2021}, Cauty \cite{cauty2005}-\cite{cauty2007} tried to solve the Schauder conjecture (see also the comment given by Dobrowolski \cite{dobrowolski2003} on Professor Cauty's work), and they (Ennassik and Taoudi \cite{ennassik2021}) gave the solution to the Schauder conjecture by a different method for single-valued continuous mappings. From the respective of development on the study of fixed point theory and related topics in nonlinear analysis, a number of works have been contributed, just mention a few of them, inclyding Agarwal et al.\cite{Agarwal}, Ben-El-Mechaiekh \cite{ben1}, Ben-El-Mechaiekh and Saidi \cite{ben2}, Browder \cite{bro1968b}, Cellina \cite{cellina}, Chang \cite{chang1997}, Ennassik et al.\cite{ennassik}, Fan \cite{fan1952}-\cite{fan1960}, G\'{o}rniewicz \cite{gorniewicz}, Granas and Dugundji \cite{granas}, Guo et al.\cite{guo2},  Nhu \cite{nhu1996}, Park \cite{park100}, Reich \cite{reich}, Smart \cite{smart}, Tychonoff \cite{tych1935},  Weber \cite{weber1992}-\cite{weber1991}, Xiao and Lu \cite{xiaolu}, Xiao and Zhu \cite{xiaozhu2011}, Xu \cite{xu2000},
Yuan \cite{yuan1998}-\cite{yuan2022c} and many other scholars, and also see the comprehensive references and related discussion under the general framework of topological vector space, or $p$-vector spaces for non-self set-valued mappings $(0 < p \leq 1)$.

The goal of this paper is to establish a general fixed point theorem for compact single-valued continuous mapping in Hausdorff $p$-vector spaces, and fixed point theorem for upper semicontinuous set-valued mappings in locally $p$-convex for $p\in (0, 1]$. These new results provide an answer to Schauder conjecture in the affirmative  under the setting of general $p$-vector spaces for compact single-valued continuous, and also give the fixed point theorems for upper semicontinuous set-valued mappings defined on $s$-convex subsets in Hausdorff locally $p$-convex spaces, which would be fundamental for nonlinear functional analysis in mathematics, where $s, p \in (0. 1]$.

The paper has three sections as follows. Section 1 is the introduction. Section 2 describes general concepts for the $p$-convex subsets of  topological vector spaces, and locally $p$-convex spaces for $p \in (0, 1]$.  In section 3, we first give a fixed point theorem for non-compact
single-valued continuous mappings in topological vector spaces, and then we establish general fixed point theorems for upper semi-continuous
mappings defined on $s$-convex subsets in locally $p$-convex compacts for $s, p \in (0, 1]$.

For the convenience of our discussion, throughout this paper, all $p$-convex vector spaces are assumed to be Hausdorff, and $p$ satisfying the the condition for  $0 <p \leq 1$ unless specified, and also we denote by $\mathbb{N}$ the set of all positive integers, i.e., $\mathbb{N}:=\{1,2, \cdots, \}$. For a set $X$, the $2^X$ denotes the family of all subsets of $X$.

\vskip.1in
\section{The basic notions and results of $p$-vector spaces}

\noindent
{\bf Definition 2.1.}
A set $A$ in a vector space $X$ is said to be
$p$-convex for $0 < p \leq 1$ if, for any $x, y\in A$, $0 \leq s, t \leq 1$ with
$s^p + t^p=1$, we have
$s^{1/p}x+t^{1/p}y\in A$; and if $A$ is $1$-convex, it is simply called convex (for $p = 1$) in general vector spaces;
the set $A$ is said to be absolutely $p$-convex if $s^{1/p}x+t^{1/p}y\in A$ for $0 \leq |s|, |t| \leq 1$ with $|s|^p + |t|^p \leq 1$.

\noindent
{\bf Definition 2.2.}
If $A$ is a subset of a topological vector space
$X$, the closure of $A$ is denoted by $\overline{A}$, then the $p$-convex
hull of $A$ and its closed $p$-convex hull are denoted by $C_{p}(A)$ and
$\overline{C}_{p}(A)$, respectively, which are the smallest $p$-convex
set containing $A$ and the smallest closed $p$-convex set containing
$A$, respectively.

\noindent
{\bf Definition 2.3.}
 Let $A$ be $p$-convex and
$x_{1}, \ldots , x_{n}\in A$, and
$t_{i}\geq 0$, $  \sum_{1}^{n}t_{i}^{\mathrm{p}}=1$. Then
$ \sum_{1}^{n}t_{i}x_{i}$ is called a $p$-convex combination
of $\{x_{i}\}$ for $i=1, 2, \ldots , n$. If
$ \sum_{1}^{n}|t_{i}|^{\mathrm{p}}\leq 1$, then
$ \sum_{1}^{n}t_{i}x_{i}$ is called an absolutely $p$-convex
combination. It is easy to see that
$ \sum_{1}^{n}t_{i}x_{i}\in A$ for a $p$-convex set
$A$.

\noindent
{\bf Definition 2.4.}
A subset $A$ of a vector space $X$ is called balanced (circled) if $\lambda A \subset A$ holds for all scalars
$\lambda $ satisfying $|\lambda | \leq 1$. By its definition, a balanced set $A$ is symmetric, i.e., $A=-A$.
We also say that the set  $A$ is absorbing
if, for each $x \in X$, there is a real number $\rho _{x} >0$ such that
$\lambda x \in A$ for all $\lambda > 0$ with
$|\lambda |\leq \rho _{x}$.

\vskip.1in
By the Definition~2.4 above, it is easy to see that the system of all circled subsets
of $X$ is easily seen to be closed under the formation of linear combinations,
arbitrary unions, and arbitrary intersections. In particular, every set
$A \subset X$ determines the smallest circled subset $\hat A$ of $X$ in
which it is contained: $\hat A$ is called the circled hull of $A$. It is
clear that $\hat A =\cup _{|\lambda |\leq 1} \lambda A$ holds, so that
$A$ is circled if and only if (in short, iff) $\hat A =A$. We use
$\overline{\hat A}$ to denote the closed circled hull of
$A\subset X$. In addition, if $X$ is a topological vector space, we use $\operatorname{int}(A)$ to denote
the interior of set $A \subset X$, and if $0 \in \operatorname{int}(A)$, then
$\operatorname{int}(A)$ is also circled, and we use $\partial A$ to denote the boundary
of $A$ in $X$ unless specified otherwise.

\noindent
{\bf Definition 2.5.} A topological vector space is said to be locally
$p$-convex if the origin has a fundamental set of absolutely $p$-convex
$0$-neighborhoods. This topology can be determined by $p$-seminorms which
are defined in the obvious way (see pp.~52 of Bayoumi \cite{bayoumi}, Jarchow
\cite{jarchow}, or Rolewicz \cite{rolewicz}). If $p=1$, $X$ is a usual locally convex space.

\vskip.1in
\noindent
{\bf Remark 2.1.} It is well-known that a~given $p$-seminorm $P$ is said to be a $p$-norm if
$x = 0$ whenever $P(x) = 0$. A~vector space with a specific $p$-norm is
called a $p$-normed space.
Specifically, a~Hausdorff topological vector space is locally bounded if and only if it is a $p$-normed space for some $p$-norm $\| \cdot \|_{p}$, where $0 < p \leq 1$ (see pp.~114 of Jarchow \cite{jarchow}).

\vskip.1in
We also note that examples of $p$-normed spaces include
$L^{p}(\mu )$ - spaces and Hardy spaces $H_{p}$, $0 < p < 1$, endowed with their usual $p$-norms. Moreover, we would like to make the following important two points:

(1) First, by the fact that (e.g., see Kalton et al.~\cite{kalton1984}) 
there is no open convex nonvoid subset in
$L^{p}[0, 1]$ (for $0< p < 1$) except $L^{p}[0,1]$ itself. This means that
$p$-normed paces with $0< p <1$ are not necessarily locally convex. Moreover,
we know that every $p$-normed space is locally $p$-convex; and incorporating
Lemma~2.2 below, 
it seems that $p$-vector spaces (for $0 < p \leq 1$) are nicer and bigger
as we can use a $p$-convex subset in locally $p$-convex spaces to approximate convex subsets in topological vector spaces (TVS) by
Lemma~2.1(ii), and also Lemma~3.1 in Section 3. In this way, it seems that $p$-vector spaces have better properties in terms of $p$-convexity than
the usually ($1-$) convex subsets used in TVS with $p=1$.

(2) Second, it is worthwhile noting that a $0$-neighborhood in a topological
vector space is always absorbing by Lemma 2.1.16 of Balachandran
\cite{balachandra} or Proposition 2.2.3 of Jarchow \cite{jarchow}.

\vskip.1in
The following result is a very important and useful result which allows
us to make the approximation for convex subsets in topological vector spaces
by $p$-convex subsets in $p$-convex vector spaces. For the reader's convenience,
we provide a sketch of proof below (see Lemma 2.1 of Ennassik and
Taoudi \cite{ennassik}, Remark 2.1 of Qiu and Rolewicz \cite{qiu}).

\vskip.1in
\noindent
{\bf Lemma 2.1.}
Let $A$ be a subset of a vector space $X$, then we have:

(i) If $A$ is $p$-convex with $0 < p < 1$, then $\alpha x \in A$ for any
$x \in A$ and any $0 < \alpha \leq 1$.

(ii) If $A$ is convex and $0 \in A$, then $A$ is $p$-convex for any
$p \in (0, 1]$.

(iii) If $A$ is $p$-convex for some $p \in (0, 1)$, then $A$ is $s$-convex
for any $s \in (0, p]$.

{\bf Proof.} See Lemma 2.1 of Ennassik and
Taoudi \cite{ennassik}, or, the Remark 2.1 of Qiu and Rolewicz \cite{qiu}). $\square$

\vskip.1in
\noindent
{\bf Remark 2.2.}
We would like to point out that results (i) and (iii)
of Lemma~2.1 do not hold for $p = 1$. Indeed, any singleton
$\{x\} \subset X$ is convex in topological vector spaces; but if
$x \neq 0$, then it is not $p$-convex for any $p \in (0, 1)$ (see also Lemma 2.2 below).

We also need the following proposition, which is Proposition 6.7.2 of Jarchow
\cite{jarchow}.

\noindent
{\bf Proposition 2.1.}
Let $K$ be compact in a topological vector
$X$ and $(1< p \leq 1)$. Then the closure $\overline{C}_{p}(K)$ of the
$p$-convex hull and the closure $\overline{AC}_{p}(K)$ of absolutely
$p$-convex hull of $K$ are compact if and only if
$\overline{C}_{p}(K)$ and $\overline{AC}_{p}(K)$ are complete, respectively.

Before we close this section, we would like to point out that the structure
of $p$-convexity when $p \in (0, 1)$ is really different from what we normally
have for the concept of ``convexity'' used in topological vector spaces.
In particular, maybe the following fact is one of the reasons for
us to use better ($p$-convex) structures in $p$-vector spaces to approximate
the corresponding structure of the convexity used in TVS (i.e., the
$p$-vector space when $p=1$).

Based on the discussion in pp.~1740 by Xiao
and Zhu \cite{xiaozhu2011}, 
we have the following fact, which indicates that each $p$-convex subset is ``bigger'' than the
convex subset in topological vector spaces for $0 < p < 1$.

\noindent
{\bf Lemma 2.2.}
Let $x$ be a point of a $p$-vector space $E$, where assume
$0 < p < 1$, then the $p$-convex hull and the closure of $\{x\}$ are given
by
%
\begin{equation}
\label{eq1} C_{p}\bigl(\{x\}\bigr)=\begin{cases}
\{tx: t \in (0, 1]\},  &\text{if }  x \neq 0,
\\
\{0\}, & \text{if }  x = 0;
\end{cases}
\end{equation}
and
%
\begin{equation}
\label{eq2} \overline{C_{p}\bigl(\{x\}\bigr)}=\begin{cases}
\{tx: t \in [0, 1]\},&  \text{if }  x \neq 0,
\\
\{0\}, & \text{if }  x= 0.
\end{cases}
\end{equation}

\vskip.1in
Here we note that if $x$ is a given one point in $p$-vector space $E$, when
$p=1$, we have that $\overline{C_{1}(\{x\})} =C_{1}(\{x\})=\{ x\}$. This
shows to be significantly different for the structure of $p$-convexity
between $p=1$ and $p\neq 1$! Thus it is necessary to establish fixed point theorems for set-valued mappings with convex values, instead of
$p$-convex values as zero is always contained by any $p$-closed subset when $p \in (0, 1)$.

Throughout this paper, without loss of generality
unless specified otherwise, for a given $p$-vector space $E$, where
$p \in (0, 1]$, we always denote by $\mathfrak{U}$ the base of the
$p$-vector space $E$'s topology structure, which is the family of its
$0$-neighborhoods, and we assume that all $p$-vector spaces $E$ are Hausdorff unless specified for $p \in (0, 1]$.
\vskip.1in

\section{Fixed point theorems in both $p$-vector spaces and locally $p$-convex spaces}
\label{sec4}

In this section, we establish a general fixed point theorem for single-valued continuous mapping in Hausdorff $p$-vector spaces, and fixed point theorem for upper semicontinuous set-valued (USC) mappings in locally $p$-convex for $p \in (0, 1]$. These new results provide an answer to Schauder conjecture in the affirmative  under the setting of general $p$-vector spaces for single-valued continuous, and upper semicontinuous set-valued mappings defined on $s$-convex subsets in locally $p$-convex spaces, which would be fundamental for nonlinear functional analysis in mathematics, where $s, p \in [0, 1]$.

Here, we first gather together necessary definitions, notations, and known facts needed in this section.

\noindent
{\bf Definition 3.1.}
Let $X$ and $Y$ be two topological spaces. A~set-valued
mapping (also called multifunction) $T: X \longrightarrow 2^{Y}$ is a point
to set function such that for each $x \in X$, $T(x)$ is a subset of
$Y$. The mapping $T$ is said to be upper semicontinuous (USC) if the subset
$T^{-1}(B): = \{ x\in X: T(x) \cap B \neq \emptyset \}$ (equivalently, the set
$\{x \in X: T(x) \subset B\}$) is closed (equivalently, open) for any closed (resp.,
open) subset $B$ in $Y$. The function $T: X \rightarrow 2^{Y}$ is said
to be lower semicontinuous (LSC) if the set $T^{-1}(A)$ is open for any
open subset $A$ in $Y$.

\vskip.1in
\noindent
{\bf Definition 3.2.}
We recall that for two given topological spaces $X$ and $Y$, and a set-valued mapping $T: X \rightarrow 2^Y$ is said to be compact if there is compact subset set $C$ in $Y$ such that $T(X)(:=\{y \in T(x), x \in X\})$ is contained in $C$, i.e., $F(X) \subset C$. Now we  have the following non-compact versions of fixed point theorems for compact single-valued mappings defined in locally $p$-convex and topological vector spaces for $ 0 < p \leq 1$.

\vskip.1in
 We now state the following result which is a compact version of Theorem 3.1 and Theorem 3.3 by Ennassik and Taoudi \cite{ennassik2021}.

\vskip.1in
\noindent
{\bf Theorem 3.1.} If $K$  is a nonempty closed $p$-convex subset of either a Hausdorff locally $p$-convex space (for $ 0 < p \leq 1$) or a Hausdorff topological vector space $X$, then the compact single-valued continuous mapping $T: K \rightarrow K$ has at least a fixed point.

\noindent
{\bf Proof.} As $T$ is compact, there exists a compact subset $A$ in $K$ such that $T(K)\subset A$.
 Let $K_0: =\overline{C}_p(A)$ be the closure of the $p$-convex hull of the subset $A$ in $K$. Then $K_0$ is compact $p$-convex by Proposition 2.1, and the mapping $T: K_0 \rightarrow K_0$ is continuous.  Then we can prove the conclusion by considering the self-mapping $T$ on $K_0$ as applications of Theorem 3.1 and Theorem 3.3 given by Ennassik and Taoudi \cite{ennassik2021}.
 $\square$

\vskip.1in
\vskip.in
Now we are going to discuss how to establish the main results for the existence of fixed point theorem for upper semicontinuous set-valued mappings defined on $s$-convex subsets  under the framework of Hausdorff locally $p$-convex spaces, where $s, p \in (0, 1]$.

\vskip.1in
By following the idea used by Repov$\breve{s}$ et al.\cite{RR} for the graph approximation of quasi upper semicontinuous set-valued mappings with the concept of the $``$$p$-convexity" used in locally $p$-convex spaces to replace the usually concept of $``$convexity" used in topological vector spaces (see also Ben-El-Mechaiekh \cite{ben1}, Ben-El-Mechaiekh and Saidi \cite{ben2}, Cellina \cite{cellina}, Kryszewski \cite{KK}, Repov$\breve{s}$ et al.\cite{RR} and related references), we have the following Lemma 3.1 which is then used as a tool to establish a general fixed point theorem for upper semicontinous set-valued mappings in Hausdorff locally $p$-convex spaces for $p \in (0, 1]$.

Here we also recall that if $X$ and $Y$ are two topological spaces and $F: X \rightarrow 2^Y$ is a set-valued mapping, and we denote by either $Graph(F)$ or $\Gamma_F$ for the graph of $F$ in $X \times Y$, and $\alpha$ is a given open cover of $\Gamma_F$ in $X \times Y$, then a (single- or set-valued) mapping $G: X \rightarrow Y$ is said to be an $\alpha$-approximation (also called $\alpha$-graph approximation) of $F$ if for each point $p \in \Gamma_G$, there exists a point $q \in \Gamma_F$ such that $p$ and $q$ lie in some common element of the over $\alpha$.
In the case $Y$ is a topological vector space, if $\Omega$ is the open cover of $X$ and $V$ is an open neighborhood of their origin in $Y$, then $\Omega \times \{y+V\}_{y \in Y}$ is one open cover of $X \times Y$, which is denoted by $\Omega \times V$ in this section. We also refer the readers the reference books by Dugungji \cite{dugundji} and Kelly \cite{kelly} for the corresponding notations and concepts used in general topology.

\noindent
{\bf Lemma 3.1.} Let $X$ be a paracompact space and $Y$ be a topological vector space and $p \in (0, 1]$. If $F: X \rightarrow 2^Y$ is an upper semicontinuous mapping with $p$-convex values, then for each open cover $\Omega$ of $X$, and each $p$-convex open neighborhood $V$ of the origin in $Y$, there exists a continuous single-valued $(\Omega \times V)$-approximation for the set-valued mapping $F$. In particular, the conclusion holds if $V$  is any convex open neighborhood of the origin in $Y$.

\noindent
{\bf Proof.} Let $\Omega$ be an open covering of $X$, and let $V$ be a $p$-convex open neighborhood of the origin in $Y$. For each $x \in X$, fix
an arbitrary element $W(x)\in \Omega$ such that $x \in W(x)$, then we first claim the following statements:

(1) By the upper semicontinuity (USC) of the mapping $F$, for each $x \in X$, there exists an open neighborhood $U(x) \subset W(x)$ such that $F(z) \subset F(x)+ V$ for all $z \in U(x)$;

(2) As $X$ is paracompact, by Theorem 3.5 of Dugundji \cite{dugundji} (see also Theorem 28 in Chapter 5 of Kelly \cite{kelly}), without loss of the generality, let the family $\{G(x)\}_{x \in X}$ be a  covering which is a star refinement of the covering $\{U(x)\}_{x \in X}$ of $X$ (and see also the discussion on pp.167-168 by Dugundji \cite{dugundji} for the concept of the star refinement for a given covering);

(3) Using the upper semicontinuity property again for the mapping $F$, for each $x \in X$, there exists an open neighborhood $U'(x) \subset G(x)$ such that $F(z) \subset F(x) + V $ for all $z \in U'(x)$;

(4) Let $\{e_{\alpha}\}_{\alpha \in A}$ be a locally finite continuous partition of unity inscribed into the covering
$\{U'(x)\}_{x \in X}$ of $X$, where $A$ is the index set, with $\Sigma_{\alpha \in A} ~e_{\alpha(x)} = 1$ for each $x \in X$; and for each $\alpha \in A$, we can choose $x_{\alpha} \in X$ such that $supp ~e_{\alpha} \subset U'(x_{\alpha})$,
 and choosing one point $y_{\alpha} \in F(x_{\alpha})$, where $supp ~e_{\alpha}$ is the support of $e_{\alpha}$ (defined by
 $supp ~e_{\alpha}:=\overline{\{x\in X: e_{\alpha}(x) \neq 0\}}$); and

(5) Finally, define a mapping $f: X \rightarrow Y$ by  $f(x):= \Sigma_{\alpha \in A} e^{\frac{1}{p}}_{\alpha}(x) y_{\alpha}$ for each $x \in X$,
where  $y_{\alpha} \in F(x_{\alpha})$ as given by (4) above,  then $f$ is well-defined, where the sum is taken over all $\alpha \in A$ with
$e_{\alpha}(x) > 0$. By (3), it follows that $\Sigma_{\alpha \in A} (e^{\frac{1}{p}}_{\alpha}(x))^p =\Sigma_{\alpha \in A} e_{\alpha}(x)=1$.

Now we show that $f$ is indeed the desired single-valued continuous mapping, which is the $(\Omega \times V)$-approximation for the mapping $F$. Indeed for any given $x_0 \in X$, we have that
$$x_0 \in St\{x_0, \{supp ~e_{\alpha} \}_{\alpha \in A} \} \subset St \{x_0, \{U'(x)\}_{x \in X}\} \subset
St \{x_0, \{G(x)\}_{x \in X}\} \subset U(x') \subset W(x')$$
\noindent
for some $x' \in X$, where $St \{x_0, \{supp ~e_{\alpha} \}_{\alpha \in A}\}$ denotes the Star of the point $\{x_0\}$ with respect to the family $ \{supp ~e_{\alpha} \}_{\alpha \in A}$ and defined by $St\{x_0, \{supp ~e_{\alpha} \}_{\alpha \in A}\}:=\cup\{U: x_0 \in U, U \in \{supp ~e_{\alpha} \}_{\alpha \in A} ~\}$ (see also the corresponding discussion for the notation and concept on pp.349 given by Ageev and Repov$\breve{s}$ \cite{ageevrepovs}).

By the definition  of quasi upper semi continuity, we have that
$x'\in W(x')$. Hence the points $x_0$ and $x'$ are $\Omega$-close.

Secondly, if $e_{\alpha}(x_0) > 0$ for $\alpha \in A$, then $x_0 \in G(x_{\alpha})$ and $x_{\alpha} \in G(x_{\alpha})$ by (3) above.
Thus $x_{\alpha}\in St\{x_0, \{G(x)\}_{x\in X} \} \subset U(x')$. Therefore,
$y_{\alpha} \in F(x_{\alpha}) \subset F(x') + V$,
i.e., $y_{\alpha} - v_{\alpha} \in V$ for some $v_{\alpha} \in F(x')$ for $\alpha \in A$.
But then, for $v: =\Sigma_{\alpha} e^{\frac{1}{p}}_{\alpha}(x_0) v_{\alpha} \in F(x')$ as $F$ is $p$-convex valued and we know that $\Sigma_{\alpha \in A} (e^{\frac{1}{p}}_{\alpha}(x))^p =\Sigma_{\alpha \in A} e_{\alpha}(x)=1$ as shown by (5) above, and $y_{\alpha} - v_{\alpha} \in V$, too, for
$\alpha\in A$, thus we have that $f(x_0)- v =\Sigma e^{\frac{1}{p}}_{\alpha}(x_0)(y_{\alpha} - v_{\alpha}) \in V$ as $V$ is $p$-convex.
Hence, the  point $(x_0, f(x_0)) \in Graph(f)$ is $(\Omega \times V)$-close to the point
$(x', v) \in Graph(F)$.

In particular, as each convex neighborhood of the origin in $Y$ is also $p$-convex for each $p \in (0, 1]$, thus the conclusion holds.  The proof is complete.  $\square$

\vskip.1in
Now we have the following main result for upper semicontinuous set-valued mappings in Hausdorff locally $p$-convex spaces.

\vskip.1in
\noindent
{\bf Theorem 3.2.} Let $K$ be a non-empty
compact $s$-convex subset of a Hausdorff locally $p$-convex space $X$, where $p, s \in (0, 1]$.
If $T: K \rightarrow 2^K$ is an upper semicontinuous set-valued mapping with non-empty closed $p$-convex values, then $T$ has a fixed point in $K$.

\noindent
{\bf Proof.} We give the proof by using the graph approximation approach for upper semicontinuous set-valued mappings established in this section above. Let $\mathfrak{U}$ be the family of absolutely $p$-convex open neighborhoods of the origin in $X$.
By the fact the family $\{x + u \}_{x \in K}$ is an open covering of $K$, and we denote the family $\{x + u \}_{x \in K}$ by $\Omega$.
Now by Lemma 3.1, it follows that there exists one (single-valued) continuous mapping $f_u: K \rightarrow K$, which is
$(\Omega \times u)$-approximation of the mapping $T$. By Theorem 3.1, $f_u$ has a fixed point $x_u = f_u(x_u)$ in $K$ for each $u \in \mathfrak{U}$. Note that $(x_u, f_u(x_u))=(x_u, x_u) \in Graph(f_u)$, which is $(\Omega \times u)$-approximation of the Graph(T), and
the graph of $T$ is closed due to the assumption, we will go to prove $T$ has a fixed point $x^*$ which is indeed the limit of
 some sub-net of the family $\{x_u\}_{u \in  \mathfrak{U}}$ in $K$, i.e., $x^* \in T(x^*)$, by using notations of language in general topology (for related references on the discussion for normed spaces or topological (vector) spaces, see Cellina \cite{cellina}, Ben-El-Mechaiekh \cite{ben1} and Fan \cite{fan1952}).

Indeed, for any given open $p$-convex member $u$ in $\mathfrak{U}$, as the set $\{x + u\}_{x \in K} \times \{y + u\}_{y \in K}$ is an open cover of $K \times K$, by Lemma 3.1, there exists a single-valued continuous mapping $f_u:  K \rightarrow K$, which is $(\Omega \times u)$-approximation of the Graph(T), where $\Omega: = \{x+ u\}_{x \in K}$ as mentioned above.
By Theorem 3.1, $f_u$ has a fixed point $x_u = f_u(x_u)$ in $K$ for each $u \in \mathfrak{U}$.
Now for $x_u \in K$, by following the proof of Lemma 3.1, we observe that firstly, there exists $x'_u \in K$ such that $x_u \in x'_u + u$;
and secondly, there also exists $ v_u \in F(x'_u)$ such that $f_u(x_u)-v_u \in u$ which means that $f_u(x_u) \in v_u + u$.

In summary, for any given $u \in \mathfrak{U}$, there exist a continuous mapping $f_u: K \rightarrow K$, which has at least one  fixed point $x_u \in K$ such that $x_u = f_u(x_u)$ with  $(x_u, x_u) =(x_u, f_u(x_u)) \in Graph(f_u)$, and we also have the following statements:

(1) there exists $x'_u \in K$ such that $x_u \in x'_u + u$; and

(2) there exist $ v_u \in F(x'_u)$ such that $f_u(x_u)-v_u \in u$, which means $f_u(x_u) \in v_u + u$.

Since $K$ is compact, without loss of the generality, we may assume that there exists a sub-net $(x_{u_i})_{u_i \in \mathfrak{U}}$ converges to $x^*$ in $K$. Now we will show that  $x^*$ is the fixed point of $T$, i.e., $x^* \in T(x^*)$.

As $K$ is compact, without loss of the generality, we may assume that two net $\{x_u\}_{u \in \mathfrak{U}}$ and $\{x'_u\}_{u \in \mathfrak{U}}$ in $K$ have the sub-net $\{x_{u_i}\}_{u_i \in \mathfrak{U}}$ converges to $x^*$, and the sub0net $\{x'_{u_i}\}_{u_i \in \mathfrak{U}}$ converges to $x'^*$, respectively in $K$. By the statement of (1) above, it is clear that we must have
 $x^* = x'^*$; otherwise, as the family $\mathfrak{U}$ is the base of absolutely $p$-convex open neighborhoods of the origin in $X$, by (1) we will have the contradiction, and thus our claim that $x^* = x'^*$ is true  in locally $p$-convex space $X$.

Now we prove that $x^*$ is a fixed point of $T$ by using the statement of (2) for all $u \in \mathfrak{U}$. As the net $\{v_u\}_{u \in \mathfrak{U}} \subset K$, we may assume its sub-net $\{v_{u_i}\}_{u_i \in \mathfrak{U}}$ converges to $v^*$. Then by the statement given by (2), it is clear that we have that $\lim_{u_i \in \mathfrak{U}}v_{u_i}=v^* = \lim_{u_i \in \mathfrak{U}}f_{u_i}(x_{u_i})=
\lim_{u_i \in \mathfrak{U}}x_{u_i} = x^*$. By the fact that $(v_{u_i}, x'_{u_i}) \in Graph(T)$, and the graph of $T$ is closed as $T$ is USC with closed values, it follows that $x^* = v^* \in T(x^*)$, which means $x^*$ is a fixed point of $T$.  The proof is complete. $\square$

\vskip.1in
\noindent
{\bf Remark 3.2.} Here we are not sure if the assumption $``T(x)$ is with non-empty closed $p$-convex values" could be replaced by the
condition $``T(x)$ is with non-empty closed $s$-convex values" in Theorem 3.2. In fact, it seems  that the proof of Theorem 4.3 given by
Ennassik et al.\cite{ennassik} only goes through for the case $s \leq p$, not for the general case when both $s, p \in (0, 1]$ (please note that the letter $p$ is denoted as the letter $r$ by Ennassik et al.\cite{ennassik}). Thus, we are still looking for a proper way to prove if the conclusion of the Theorem 3.2 is true under Hausdorff topological vector spaces instead of locally $p$-convex spaces.

\vskip.in
As an immediate consequence of Theorem 3.2, we have following fixed point result for upper semicontinuous set-valued mappings in locally convex spaces for compact $s$-convex subsets, which include the common compact convex sets as a special class when $s=1$.

\vskip.1in
By following the same idea used in the proof of Theorem 3.1, we have the following Theorem 3.3 for compact upper semicontinuous set-valued mappings and thus we omit its proof here.

\vskip.1in
\noindent
{\bf Theorem 3.3.} If $K$ is a nonempty closed $s$-convex subset of a Hausdorff locally $p$-convex space $X$, where $s, p \in (0, 1]$, then any compact upper semicontinuous set-valued mapping $T: K \rightarrow 2^K$  with nonempty closed $p$-convex values, has at least one fixed point.
\vskip.1in

\vskip.1in
\noindent
{\bf Corollary 3.1.} If $K$ is a nonempty closed $s$-convex subset of a Hausdorff locally convex space $X$, where $s \in (0, 1]$, then any upper semicontinuous set-valued mapping $T: K \rightarrow 2^K$ with non-empty closed convex values, has at least one fixed point.

\noindent
{\bf Proof.} Let $p=1$ in Theorem 3.3, then the conclusion follows by Theorem 3.3. This completes the proof.
$\square$.

\vskip.1in
As a special case of Theorem 3.3 or Corollary 3.1, we also have the following corollary.

\vskip.1in
\noindent
{\bf Corollary 3.2.} If $K$ is a nonempty compact $s$-convex subset of a Hausdorff locally convex space $X$, where $s \in (0, 1]$, then any upper semi-continuous set-valued mapping $T: K \rightarrow 2^K$ with non-empty closed convex values, has at least one fixed point.

\vskip.1in
\noindent
{\bf Remark 3.3.}
Theorem 3.1 says that each compact single-valued mapping defined on a closed $p$-convex subsets ($0< p \leq 1)$ in topological vector spaces has the fixed point property, which does not only include or improve most available results for fixed point theorems  in the existing literature as special cases (just mention a few, Ben-El-Mechaiekh \cite{ben1}, Ben-El-Mechaiekh and Saidi \cite{ben2}, Ennassik and Taoudi \cite{ennassik2021}, Mauldin \cite{Mauldin}, 
Granas and Dugundji \cite{granas}, O'Regan and Precup \cite{OP}, Reich \cite{reich}, Park \cite{park100} and references wherein). In particular, we note that the answer to Schauder conjecture in the affirmative for a single-valued continuous mapping recently was obtained by Ennaassik and Taoudi \cite{ennassik2021} defined on non-empty compact $p$-convex subset in Hausdorff topologocal vector spaces, where $p \in (0, 1]$.


We first note that Theorem 3.3 improve or unifies corresponding results given by
Cauty \cite{cauty2005}, Cauty \cite{cauty2007}, 
Dobrowolski \cite{dobrowolski2003}, Nhu \cite{nhu1996}, Park \cite{park100}, Reich \cite{reich}, Smart \cite{smart}, Xiao and Lu \cite{xiaolu}, Xiao and Zhu \cite{xiaozhu2011}, Yuan \cite{yuan1999}-\cite{yuan2022c} under the framework of compact single-valued or upper semicontinuous set-valued mappings.

We also like to mention that by comparing with topological degree approach or related other methods used or developed by 
Cauty \cite{cauty2005}-\cite{cauty2007}, Nhu \cite{nhu1996} and others, the arguments used in this section actually provides an accessible way for the study of nonlinear analysis for $p$-convex vector spaces $(0 < p \leq 1$). The results given in this paper are new, and may easily understand for general readers in mathematical community, and see more by Yuan \cite{yuan2022b}-\cite{yuan2022c} and related references on the study of nonlinear analysis and related applications in both $p$-vector and locally $p$-convex spaces for $0 < p \leq 1$.

{\bf Acknowledgement}


This research is partially supported by the National Natural Science Foundation of China [grant numbers 71971031 and
U1811462].

\vskip.1in
\noindent
{\bf Compliance with Ethical Standards}

The author declares that there is no conflict of interest.

\end{document}